\input amstex
\pageno=1

\magnification=1200
\loadmsbm
\loadmsam
\loadeufm
\input amssym
\UseAMSsymbols

\TagsOnRight

\hsize172 true mm
\vsize212 true mm
\voffset=20 true mm
\hoffset=0 true mm
\baselineskip 7.0 true mm plus0.4 true mm minus0.4 true  mm
\def\st{\rightarrow}
\centerline{\bf GEOMETRY AND ALGEBRA OF REAL FORMS OF COMPLEX CURVES}
\bigskip
\centerline{\bf S.M.NATANZON}
\bigskip
\centerline{\bf Moscow State University and
Independent University of Moscow}
\vskip 1cm

\noindent {\bf Introduction}
\vskip 0.4cm
\noindent
Let $X$ be a non-singular, non-reducible real algebraic curve of genus
$g=g(X)$. It is called orientable if its real points $\Bbb R (X)\ $ divide its
complex points $\Bbb C(X)\supset \Bbb R (X)$ into 2 connected components.
Consider now a set of
orientable non-singular,  non-reducible real algebraic curves
$X_1,\dots, X_n$ $(n>3)$ of genus $g>1$ such  that  for  all  $i\ne  j$
$X_i$ is non-isomorphic to $X_j$ over $\Bbb R$ but it is isomorphic to $X_j$
over $\Bbb C$. According to [3], $(n-4)2^{n-3}\leqslant g-1$ and there exists
$\{X_1,\dots, X_n\}$ such that $(n-4)2^{n-3}=g-1$.
According to Harnak theorem,  $\Bbb R(X_i)$
form $\vert X_i\vert\leqslant g+1$ simple closed contours (ovals).

In this paper we prove that
     $$\sum^n_{i=1}\vert X_i\vert\leqslant 2g-(n-9)2^{n-3}-2\leqslant
2g+30$$
and these estimates are exact.

For $n=3$ and 4 it was proved in [4].

The proof is based on a detal description of
real forms of complex algebraic curves. This
description has self-dependent importance .
By our conditions all complex
algebraic curves $P_i=\Bbb C(X_i)$ are isomorphic to a complex algebraic
curve $P$ of genus $g$. Consider biholomorphic maps
$\varphi_i:P_i\rightarrow P$. The involutions of complex conjugations
$\tau'_i:P_i\rightarrow P_i$  give antiholomorphic involutions
     $\tau_i=\varphi\tau'_i\varphi^{-1} :P\rightarrow P$.
They generate a finite group $W:P\rightarrow P$. In $\S$1
following [5] we prove that $W$ is a Coxeter group.
In $\S$2 we give a complete description of
all such pairs $(P,W)$. In $\S$3, using the results of $\S$1, $\S$2,
and the classification of finite Coxeter groups we prove that
     $$\sum^n_{i=1}\vert X_i\vert \leqslant 2g-(n-9)2^{n-3}-2.$$
Farther, using an example of Singerman [8], we construct
families of orientable curves $\{X_i,\dots,$ $X_n\}$ such that
     $$\sum^n_{i=1}\vert X_i\vert= 2g-(n-9)2^{n-3}-2.$$

Some of these results were announced in [6, 7].

This work was carried out with the financial support by grants: 
RFBR 98-01-00612 and INTAS 96-0713.

\bigskip
\noindent
{\bf 1. Real equipments and Coxeter groups}
\bigskip
\noindent
Let $P$ be a complex algebraic curve, that is a compact Riemann surface
of genus $g(P)$. An antiholomorphic involution $\tau:P\st P$ is called
{\sl a real form} of $P$ $[4]$. It gives {\sl a real algebraic curve} $(P,\tau)$
$[1]$ with {\sl real points}
$$P^\tau=\{p\in P\vert \tau p=p\}.$$ It is obvious that if the set
$P^{\tau_1}\cap P^{\tau_2}$ is infinite then $\tau_1=\tau_2$.

A real form $\tau$ is called {\sl orientable} if $P/<\tau>$ is an orientable
surface. In this case $P^\tau$ divides $P$ into 2 connected components.

A finite group $W$ generated by real (orientable) forms is called
{\sl a real (orientable) equipment} of $P$.

Let $W$ be a real orientable equipment of $P$. Denote by
$[W]$ the set of all real orientable forms $\tau\in W$.
The closure $C$ of a connected component of
$P\setminus \cup_{\tau\in [W]}P^\tau $
is called a {\sl camera} of $W$. Let us consider some camera $C$.
The {\sl basis} of $C$ is the set $S_C$ of all $\sigma\in [W]$ such that
the set $P^\sigma\cap C$ is infinite.
\vskip 0.4cm
\noindent
{\bf Lemma 1.1.} {\sl The basis $S_C$ generates $W$. Each $\tau\in [W]$ is
conjugated to some $\sigma\in S_C$.}
\vskip 0.2cm
\noindent{\sl Proof:}Let $\tau\in [W]$.
 Let us consider the group $W_C$ generated by $S_C$. The set
$$\widetilde P=\bigcup_{\widetilde w\in W_C} \widetilde w C\subset P$$ is
 compact. It has no boundary and therefore $\widetilde P=P$.
Thus there exists $\widetilde w\in W_C$ such that the set
$\widetilde w C\cap P^\tau$ is infinite. Furthermore
$$\widetilde w C\cap P^\tau\subset \partial (\widetilde w C)=\widetilde w
(\partial C)\subset \widetilde w(\bigcup_{\sigma\in S_C} P^\sigma).$$
Therefore there exists $\sigma\in S_C$ such that the set
$P^\tau\cap \widetilde w(P^\sigma)$ is infinite. Thus the set
$$P^\tau\cap P^{\widetilde w\sigma\widetilde w^{-1}}$$ is infinite
and $\tau=\widetilde w \sigma \widetilde w^{-1}$. It follows that $W_C$
contains $[W]$ and therefore $S_C$ generates $W$. $\square$

Let $\sigma\in [W]$ and
$W_\sigma$ be the set of all $w\in W$ such that $C$
and $w(C)$ belong to
the same connected component of $P\setminus P^\sigma$.
\vskip 0.4cm
\noindent {\bf Lemma 1.2.} {\sl
 Let $\sigma_1, \sigma_2\in S_C$, $w\in W_{\sigma_1}$ and
$w\sigma_2\notin W_{\sigma_1}$. Then $w\sigma_2=\sigma_1 w$.}
\vskip 0.2cm
\noindent{\sl Proof:}
The sets $w\sigma_2(C)$ and $C$ belong to opposite connected components of
$P\setminus P^{\sigma_1}$. The sets $C$ and $w(C)$ belong to the same connected
component of $P\setminus P^{\sigma_1}$. Thus $w\sigma_2(C)$ and $w(C)$
belong to the opposite connected components of $P\setminus P^{\sigma_1}$.
Therefore $\sigma_2(C)$ and $C$ belong to the opposite connected components of
$P\setminus P^{w^{-1}\sigma_1 w}$ and hence
$$\sigma_2(C)\cap C\subset P^{w^{-1}\sigma_1 w}.$$
On the other hand
$$\sigma_2(C)\cap P^{\sigma_2}= C\cap P^{\sigma_2}$$ and thus the set
$$\sigma_2(C)\cap C=\sigma_2(C)\cap C\cap P^{\sigma_2}=C\cap P^{\sigma_2}$$
is infinite. Therefore the set
$$P^{\sigma_2}\cap P^{w^{-1}\sigma_1 w}\supset \sigma_2(C)\cap C$$
is infinite  and $\sigma_2=w^{-1}\sigma_1 w$. $\square$

Let  $l(w)$ be the least $l$ such that $w =\sigma_1\dots \sigma_l$, where
$\sigma_i\in S_C$.
\vskip 0.4cm
\noindent {\bf Theorem 1.1} [5]. {\sl Let $W$ be a real orientable
equipment of $P$ and $C$ its camera.
Then: 1)Pair $(W,S_C)$ is a Coxeter system that is to say
$S_C=\{\sigma_1,\dots, \sigma_n\}$ generates $W$ with defining relations
$\sigma^2_i=1$, $(\sigma_i\sigma_j)^{m_{ij}}=1$ for some integers
$m_{ij}$; 2)$C$ is a fundamental region of $W$.}
\vskip 0.2cm
\noindent
{\sl Proof:} It is obvious that $1\in W_\sigma$ and
$W_\sigma\cap \sigma W_\sigma=\emptyset$ for $\sigma\in S_C$.
It follows from
$[2, IY, \S 1, n^o 7]$ that these relations and the proposition of
lemma 1.2 give that $(W, S_C)$ is a Coxeter system and
$$W_\sigma=\{w\in W\vert l(\sigma w)>l(w)\}.$$ Thus if $w\in W$ and
$wC=C$ then
$$w\in \bigcap_{\sigma\in S_C} W_\sigma=\{w\in W\vert l(\sigma w)>l(w)
\quad\text{for all}\quad \sigma\in S_C\}=1.$$
Since $C\cap w(\partial C)\subset\partial C$, we see that $C$
is a fundamental region of $W$.$\square$
\vskip 0.4cm
\noindent
{\bf Corollary 1.1.} {\sl Let $W$ be a real orientable equipment of $P$ and
$C$ its camera. Then: 1) All fixed points of all $w\in W$ belong to
$\bigcup_{\tau\in [W]} P^\tau$;
 2) If a real form $\tau\in W$ has real points then  it is conjugated
to some  $\sigma\in  S_C$;  3) If order of $w\in W$ is more that 5 and
$w$ has a fixed point, then $w$ generates a normal subgroup of $W$. }
\vskip 0.2cm
\noindent
{\sl Proof:}
1) Suppose $w\in W$ has a fixed point
$p\in P\setminus \bigcup_{\tau\in [W]} P^\tau$.
 Since $C$ is a fundamental region of $W$
and $$W(\bigcup_{\tau\in [W]} P^\tau)=\bigcup_{\tau\in [W]} P^\tau,$$
there exists
$h\in W$ such that $hwh^{-1}$ has a fixed point in $C\setminus\partial C$.
Hence $hwh^{-1}(C)=C$. Therefore $hwh^{-1}=1$ and $w=1$.
2) Let $\tau\in W$ be a real form with real points. 
Then, by 1),  $\tau\in [W]$.
It follows from this and lemma 1.1 that $\tau$ is conjugated to some
$\sigma\in S_C$.  3)  Let  $wp=p$.  Consider  a  camera  $C\ni p$. Let
$S_C=\{\sigma_1,\dots,\sigma_n\}$. It follows from 1) that
$p\in P^{\sigma_i}\cap P^{\sigma_j}$. Thus $w=(\sigma_i\sigma_j)^k$ and
order of $\sigma_i\sigma_j$ is greater that 5.  It follows from this and
the classification of
Coxeter systems  [2, VI,  $\S$4] that $\sigma_i\sigma_j$
generates a normal subgroup of $W$.  Therefore $w$ generates a  normal
subgroup of $W$.

\bigskip
\noindent
{\bf 2. Topological classification and uniformization of equipments}
\bigskip
\noindent
Let $(W,S)$ be a Coxeter system and $m_1,\dots,m_k$ be positive
integer numbers. We shall say that a set $(W,S,T)$
is $(m_1,\dots,m_k)$ - {\sl swelling Coxeter system} if
$$T=\{\sigma(i,j)\in S\vert i=1,\dots,k, j\in Z\},$$ where

$$\bigcup_{ij}\sigma(i,j)=S, \quad
\sigma(i,j+m_i)=\sigma(i,j), \ \text{and} \quad
\sigma(i,j)\ne\sigma(i,j+1) \ \text{if}\ m_i>1.$$

We say that a ($m_1^1,\dots,m^1_k)$ - swelling Coxeter
system $(W^1,S^1,T^1)$ is
{\sl isomorphic} to a $(m^2_1, \dots, m^2_k)$ - swelling Coxeter system
$(W^2,S^2,T^2)$ if there exists a permutation
$\eta:\{1,\dots,k\}\st\{1,\dots,k\}$, integers $t_1,\dots,t_k\in Z$,
and an isomorphism $\psi:W^1\st W^2$ such that $m^1_i=m^2_{\eta(i)}$ and
$$\psi(\sigma^1(i,j))=\sigma^2(\eta(i),j+t_i),$$
where $T^l=\{\sigma^l(i,j)\}$.

Let us now associate with every real orientable equipment $(P,W)$
some swelling Coxeter system.

Let $C\subset P$ be a camera of $W$ and $a_1,\dots,a_k$ be connected components
of $\partial C$. The complex structure on $P$ gives an orientation on $C$.
This orientation on $C$ gives the
orientations on $a_i$. 

If $a_i$ is not oval thus points of intersections of ovals divide 
$a_i$ on segments $l_i^1,\dots, l_i^{m_i}$.
We label the segments by the index $j$ in order
that $l^j_i\cap l^{j+1}_i \ne \emptyset$, and the
ordering of the segments $l_i^1,l_i^2,\dots,l_i^{m_i}$
gives the contour $a_i$ with the given orientation on it.

If $a_i$ is an oval that put us $m_i=1$, $l^1_i=a_i$.
For any $l^j_i$ it exists $\sigma(i,j)\in S_C$ such
that $$l_i^j\subset P^{\sigma(i,j)}.$$
Moreover $\sigma(i,j)\ne\sigma(i,j+1)$ and
$\sigma(i,1)\ne\sigma(i,m_i)$ for $m_i\ne 1$.
Put
$$T_C=\{\sigma(i,j)\vert i=1,\dots,k, \quad j\in Z\},$$ where
$\sigma(i,j+nm_i)=\sigma(i,j)$ for $n\in Z$. It is obvious that $(W, S_C, T_C)$
is a $(m_1,\dots, m_k)$ - swelling Coxeter system.

Two real orientable equipments $(P^1,W^1)$ and $(P^2,W^2)$ are called
{\sl topological equivalent} if there exists a homeomorphism
$\varphi:P^2\st P^1$ such that $W^2=\varphi W^1\varphi^{-1}$.

\vskip 0.4cm
\noindent {\bf Theorem 2.1.} {\sl
 Real orientable equipments $(P^1,W^1)$ and
$(P^2,W^2)$ are topologically equiva\-lent if and only if
$$g(P^1/W^1)=g(P^2/W^2)$$ and there exist cameras $C^l\subset P^l$
of $W^l$ such that the swelling Coxeter systems $$(W^1, S_{C^1}, T_{C^1});
\quad (W^2, S_{C^2}, T_{C^2})$$ are isomorphic.}
\vskip 0.2cm
\noindent
{\sl Proof:} Let  $(P^1,W^1)$ and $(P^2,W^2)$ be
topologically equivalent and $\varphi :P^1\st P^2$ be a homeomorphism
such that $W^2=\varphi W^1\varphi^{-1}$. Let $C^1\subset P^1$ be
a camera of $W^1$ and $C^2=\varphi(C^1)$. Consider a homomorphism
$\psi:W^1\st W^2$
such that $\psi(w)=\varphi w \varphi^{-1}$. Then it is obvious,
that $g(P^1/W^1)=g(P^2/W^2)$ and $\psi$ gives an isomorphism between
$(W^1,S_{C^1}, T_{C^1})$ and $(W^2,S_{C^2}, T_{C^2})$.

Let us now suppose $g(P^1/W^1)=g(P^2/W^2)$
and $(W^1, S_{C^1}, T_{C^1})$, $(W^2, S_{C^2},$ $T_{C^2})$ be isomor\-phic
swelling Coxeter systems for some cameras $C^l\subset P^l$.
The boundaries of $C^l$ consist of segments
$$l^{l}_{ij}\subset P^{\sigma^l(i,j)},$$ where $T_{C^l}=\{\sigma^l(i,j)\}$.
The isomorphism
$$\psi:(W^1, S_{C^1}, T_{C^1})\st (W^2, S_{C^2}, T_{C^2})$$
gives a correspondence $(i,j)\mapsto (\eta(i),\xi(j))$ such that
$$\psi(\sigma^1(i,j))=\sigma^2(\eta(i), \xi(j))$$ and
$$\sigma^2(\eta(i), \xi(j+1))=\sigma^2(\eta(i),\xi(j)+1).$$
Thus there exists a homeomorphism $\widetilde\varphi:C^1\st C^2$ such that
$$\widetilde\varphi(l^1_{ij})= l^2_{\eta(i)\xi(j)}.$$
Consider now the homeomorphism
$\varphi:P^1\st P^2$, where
$$\varphi(p)=\psi(w)\widetilde\varphi(w^{-1}p)$$ for $p\in wC^1$.
Then $W^2=\varphi W^1\varphi^{-1}$. $\square$

Let $(W,S,T)$ be a $(m_1,\dots,m_k)$ - swelling Coxeter system and
$T=\{(\sigma(i,j)\}$. Let $n(i,j)$ be order of
$(\sigma(i,j)\cdot\sigma(i,j+1))$. Put
$$\mu_g=\mu_g(W,S,T)=4g+2k-4+\sum_{i=1}^k\sum_{j=1}^{m_i}(1-\frac{1}{n(i,j)}).$$
Denote
$$\Lambda=\Lambda((W,S,T),g)$$  the Riemann sphere  if $\mu_g=0$,
the complex plane $\Bbb C$ if $\mu_g=1, $ and the upper half-plane
$$\{z\in \Bbb C\vert \text{Im} \ z>0\} $$ if $\mu>1$. Let
$\overline{\text{Aut}}(\Lambda)$
be the group of all holomorphic and all antiholomorphic
automorphisms of $\Lambda$ and
$$\text{Aut}(\Lambda)\subset \overline{\text{Aut}}(\Lambda)$$
be the subgroup of holomorphic automorphisms.

Consider a discrete group $G\subset\overline{\text{Aut}}(\Lambda)$
and an epimorphism $\psi:G\st W$. The pair $(G,\psi)$
is called a {\sl $g$-planar realization of $(W,S,T)$} if $G$ has generators
$$\{a_\alpha,b_\alpha\in \text{Aut}(\Lambda)\quad (\alpha=1,\dots,g),\quad
c_i\in \text{Aut}(\Lambda)\quad (i=1,\dots,k), $$
$$\sigma_{ij}\notin\text{Aut}(\Lambda) (i=1,\dots,k, j=1,\dots,m_i+1)\},$$
generating $G$ with defining relations
$$\prod_{\alpha=1}^g[a_\alpha b_\alpha]\prod_{i=1}^k c_i=1, \quad
\sigma^2_{ij}=1, \quad (\sigma_{ij}\cdot\sigma_{i j+1})^{n(i,j)}=1,\quad
\sigma_{i1} c_i \sigma_{im_i+1}=c_i$$
and moreover
$$\psi(\sigma_{ij})=\sigma(i,j), \quad \psi(a_\alpha)=\psi(b_\alpha)=
\psi(c_i)=1.$$ Here $[ab]=aba^{-1} b^{-1}$ is the commutator of
the elements $a,b\in G$.

\vskip 0.4cm
\noindent {\bf Lemma 2.1.} {\sl
 Let $(G,\psi)$ be a $g$-planar realization of
$(W,S,T)$, $P=\Lambda/\text{ker}\ \psi$. $W_P=G/\text{ker}\ \psi$.
Then $(P,W_P)$ is a real orientable equipment and there exists a camera
$C\subset P$ such that the swelling Coxeter system $(W_P,S_C,T_C)$
is isomorphic to $(W,S,T)$.}

\vskip 0.2cm
\noindent
{\sl Proof:}
Put $\widetilde P=\Lambda/G$. Let $F:G\st W_{P}$, $\Phi:\Lambda\st P$,
$\varphi:P\st\widetilde P$ and
$\tilde\Phi=\varphi\Phi:\Lambda\st \widetilde P$
be the natural projections. It  follows from
[9, Ch4] that all critical values of
$\widetilde\Phi$ belong to $\partial\widetilde P$.
Moreover the fundamental group of $\widetilde P$ is generated by
the images of
$a_\alpha, b_\alpha, c_i$. Thus $\varphi$ is a homeomorphism on each
connected component of
$\varphi^{-1}(\widetilde P\setminus\partial\widetilde P)$
and the closing of each one is a fundamental region of $W_P$.
It follows from  [9, Ch4] that $G$ has a fundamental region $B$
such that $$\partial\widetilde\Phi(B)= \widetilde\Phi
(\bigcup_{ij}\{z\in B\vert \sigma_{ij} z=z\}).$$
Thus the  $C=\Phi(B)$ is a fundamental
region of $W_P$ and
$$\partial C\subset\bigcup_{ij} P^{F(\sigma_{ij})}.$$
Put $$S^1_C=\{\sigma_1,\dots,\sigma_n\}=F(\{\sigma_{ij}\
(i=1,\dots,k; \ j=1,\dots,m_i\}).$$
Then $(W_P,S^1_C)$ is a Coxeter system and
$$\partial C\subset \bigcup^m_{i=1} P^{\sigma_i}.$$

Let us prove that each real form $\sigma_s$ is orientable.
Consider $$W_s=\{w\in W_P\vert l(\sigma_sw)>l(w)\},$$
where $l(w)$ is the least $l$ such that
$w=\sigma_{i_1}\dotsb\sigma_{i_l}$.
Put $$P_1=W_s C\quad \text{and} \quad P_2=(W\setminus W_s)C.$$
Let $$p\in P_1\cap P_2=\partial P_1=\partial P_2.$$
Then $p\in w_1C\cap w_2 C$, where $w_1\in W_s$,
$w_2\in W\setminus W_s$. Put $p_0=w^{-1}_1p\in \partial C$.
Then $\sigma_tp_0=p_0$ for some $\sigma_t\in S^1_C$.
Thus $$(w_1\sigma_tw_1^{-1})p=p \quad\text{and} \quad
(w_1\sigma_tw^{-1}_1)(w_1C)=w_2C.$$
Therefore $w_1\sigma_t=w_2\notin W_s$ and according to
[2, IY, $\S 1 \quad n^0 7$] $w_1\sigma_tw^{-1}_1=\sigma_s$.
Thus $$p\in P^{\sigma_s}\quad \text{and}\quad P_1\cap P_2\subset P^{\sigma_s}.$$
Therefore $P^{\sigma_s}$ divides $P$ into 2 connected components.

Thus $(P,W_P)$ is a real orientable equipment with the camera
$C$,  $S^1_C=S_C$ and $\psi$ gives an isomorphism between
$(W_P,S_C,T_C)$ and $(W,S,T)$.$\square$
\vskip 0.4cm
\noindent {\bf Theorem 2.2.} {\sl  For each $(m_1,\dots,m_k)$-swelling
Coxeter group $(W, S, T)$ and $g\ge 0$ there exists
a $g$ - planar realization of $(W,S,T)$ and a real
orientable equipment $(P,W_P)$ with a camera $C\subset P$ such
that $g(P/W_P)=g$ and the swelling Coxeter system $(W_P,S_C,T_C)$
is isomorphic to $(W,S,T)$.}

\vskip 0.2cm
\noindent
{\sl Proof:} It follows from [9, Ch 4] that there exists a discrete group
$$G\in \overline{\text{Aut}}(\Lambda)\quad (\Lambda=\Lambda((W,S,T),g)$$ with
generators
$$\{a_\alpha,b_\alpha\in \text{Aut}(\Lambda)\quad(\alpha=1,\dots,g),
\qquad c_i\in \text{Aut}(\Lambda)\quad(i=1,\dots,k),$$
$$\sigma_{ij}\notin \text{Aut}(\Lambda)\quad(i=1,\dots,k,\quad j=1,\dots,
m_i+1)\}$$
and the defining relations
$$\sigma^2_{ij}=1, \quad (\sigma_{ij}\cdot
\sigma_{ij+1})^{n(i,j)}=1,$$
$$ \prod_{\alpha=1}^g[a_\alpha b_\alpha]
\prod^k_{i=1} c_i=1, \quad \sigma_{i1} c_i\sigma_{im_i+1}=c_i,$$
where $n(i,j)$ is the order of
$$(\sigma(i,j)\cdot \sigma(i,j+1))$$ and $$T=\{\sigma(i,j)\}.$$
Put now $$\psi(a_\alpha)=\psi(b_\alpha)=\psi(c_i)=1$$ and
$$\psi(\sigma_{ij})=\sigma(i,j).$$ Then $(G,\psi)$ is a $g$-
planar realization
of $(W,S,T)$ and theorem 2.2 follows from lemma 2.1. $\square$

We say that real equipments $(P_1,W_1)$ and $(P_2,W_2)$ are
{\sl isomorphic}
if there exists a holo\-morphic map $\varphi:P_1\st P_2$ such that
$W_2=\varphi W_1\varphi^{-1}$.
\vskip 0,4cm
\noindent {\bf Theorem 2.3.} {\sl Each  real orientable equipment is
isomorphic to
$$(\Lambda/ \text{ker}\ \psi, G/\text{ker} \ \psi),$$ where
$(G,\psi)$ is a $g$-planar realization of some swelling Coxeter system.}
\vskip 0.2cm
\noindent
{\sl Proof:} Let $(P,W_P)$ be a real orientable equipment, $g=g(P/W_P)$
and $C$
be some camera of $W_P$.
It follows from theorem 2.2 that there exist a
 $g$- planar realization $(G_0,\psi_0)$
of the swelling Coxeter system
$(W_P, S_C, T_C)$.
It follows from lemma 2.1 and theorem 2.1 that there exists a
homeomorphism $$\varphi:P\st\Lambda/\text{ker}\  \psi_0$$ such that
$$\varphi W_P \varphi^{-1}=G_0/\text{ker}\ \psi_0.$$ Consider a
uniformization $\Phi:\Lambda\st P$ and a natural projection
$$\Phi_0:\Lambda\st\Lambda/\text{ker}\ \psi_0.$$
According to [9, Ch.5], there exists a homeomorphism
$\widetilde\varphi:\Lambda\st \Lambda$ such that
$\Phi_0\widetilde\varphi=\varphi\Phi$. Put
$G=\widetilde\varphi^{-1}G_0\widetilde\varphi$ and $\psi=\psi_0 F$,
where $F:G\st G_0$ and $F(w)=\widetilde\varphi w \widetilde\varphi^{-1}$.
Then $(G,\psi)$ is a $g$-realization of $(W_P,S_C,T_S)$ and $\Phi$
gives a isomorphism between
$(\Lambda/\text{ker}\ \psi,\ G/\text{ker}\ \psi)$
and $(P,W_P)$. $\square$
\vskip 0,4cm
\noindent
{\bf 3.  Total  number  of ovals}
\vskip 0,4cm
\noindent
We shall use
\vskip 0,4cm
\noindent {\bf Lemma 3.1.}[4] {\sl Let $W:P\rightarrow P$ be a  finite
group of autohomeomorphisms of a compact orientable surface $P$.  Then
there exists a complex structure on $P$ such that $W$ consists of
holomorphic and antiholomorphic homeomorphisms.}

It follows from Harnak's theorem that the set  $P^\alpha$ of the fixed points
of a real form $\alpha:P\rightarrow P$ consists of  $\vert\alpha\vert\leqslant
g(P)+1$ simple  closed contours (ovals).  For a real orientable equipment
$(P,W)$
we put $$h(P,W)=\sum_{\alpha\in [W]} \vert\alpha\vert -2g(P)$$
and $$P^W=\bigcup_{\alpha\in [W]} P^\alpha.$$

A real  orientable equipment $(P,W)$ is called {\sl commutative} if $W$ is a
commutative group. Put $h(n,g)=\text{max} \,\{h(P,W)\}$, where $\max$ is
taken by all commutative equipments $(P,W)$ such that $g(P)=g$ and $[W]$
consists of $n$ elements.

Consider a function $f(2)=2$,  $f(n)=-(n-9)2^{n-3}-2$ for  $n>2$.  Our
next goal is the proof that $h(n,g)\leqslant f(n)$ for $n>1$.

\vskip 0.4cm
\noindent
{\bf Lemma  3.2.}  {\sl  Suppose $h(n',g')\leqslant f(n')$ for all $g'<g$
and $n'>2$. Let $(P,W)$ be a commutative equipment such that
$g(P)=g$, $[W]=(\alpha_1,\dots,\alpha_n)$, where $n>2$, and some
connected component of $P^W$ belongs to $P^{\alpha_n}$. Then
$h(P,W)<f(n)$.}

\noindent
{\sl Proof:} Let $a$ be a connected component of $P^W$ and
$a\subset P^{\alpha_n}$. Then $a$ is an oval of $\alpha_n$ and
$A=\bigcup_{w\in W}  w(a)$ consists of $2^{n-1}$  non-intersecting
contours. If $P\setminus A$ is connected, we squeeze to
a point each boundary contour of $P\setminus A$ . Thus we obtain a compact
surface $P'$,   where   $g'=g(P')=g-2^{n-1}$.   The    forms
$\alpha_1,\dots,\alpha_n$ give   involutions
$\alpha'_i:P'\rightarrow P'$. It follows from lemma 3.1 that it exists
a complex structure on $P'$ such that $\alpha'_1,\dots,\alpha'_n$ are
real forms of $P'$. They generate a commutative equipment such that
$$\sum^n_{i=1}\vert\alpha'_i\vert=\sum^n_{i=1}
\vert\alpha_i\vert-2^{n-1}.$$ Therefore
$$h(P,W)=\sum^n_{i=1}\vert\alpha_i\vert-2g=\sum^n_{i=1}
\vert\alpha'_i\vert+2^{n-1}-2g'-2^n=(\sum^n_{i=1}\vert\alpha'_i\vert
-2g')-2^{n-1}<h(n,g')\leqslant f(n).$$

Let now $P\setminus A$ be disconnected. Then $P\setminus A$
consists of two connected components $P_1$ and $P_2$. Contract to
a point each boundary contour of $P_1$ to produce a compact
surface $P'$ of genus $$g'=g(P')=\frac{1}{2}(g-2^{n-1}+1).$$
The forms   $\alpha_1,\dots,\alpha_{n-1}$   give   involutions
$\alpha'_i:P'\rightarrow P'$.  Consider a complex  structure  on  $P'$
such that $\alpha'_1,\dots,\alpha'_{n-1}$ are real forms of $P'$. They
generate a commutative equipment such that
     $$\sum^n_{i=1}\vert\alpha'_i\vert=\frac{1}{2}(\sum^n_{i=1}
\vert\alpha_i\vert-2^{n-1}).$$
Therefore $$h(P,W)=\sum_{i=1}^n\vert\alpha_i\vert-2g=2\sum^{n-1}_{i=1}
\vert\alpha'_i\vert+2^{n-1}-4g'-2^n+2=$$
$$2(\sum^{n-1}_{i=1}\vert\alpha'_i
\vert-2g')-2^{n-1}+2\leqslant 2f(n-1)-2^{n-1}+2.$$
     Thus if $n=3$,  then $$h(P,W)\leqslant 4-4+2<4=f(3).$$ If  $n>3$,
then $$h(P,W)\leqslant
2(-(n-10)2^{n-4}-2)-2^{n-1}+2=-2^{n-3}(n-6)-2<f(n).\square$$

\vskip 0.4cm
\noindent
{\bf Lemma 3.3.} {\sl Suppose $h(n',g')\leqslant f(n')$ for  all $g'<g$ and
 $n'>2$.  Let $(P,W)$ be a commutative equipment such that $g(P)=g$,
$[W]=(\alpha_1,\dots,\alpha_n)$, where $n>2$, and some connected
component of $P^W$ belongs to $P^{\alpha_{n-1}}\cup P^{\alpha_n}$. Then
$f(P,W)<f(n)$.}

\noindent
{\sl Proof:}  Let  $a$ be a connected component of $P^W$ and $a\subset
P^{\alpha_{n-1}}\cup P^{\alpha_n}$.  If $a\subset P^{\alpha_{n-1}}$  or
$a\subset P^{\alpha_n}$,  then  lemma 3.3 follows from lemma 3.2.  Let
$a\not\subset P^{\alpha_{n-1}}$ and  $a\not\subset  P^{\alpha_{n}}$.  In
this case $a$ consists of some number $m$ of ovals from
$ P^{\alpha_{n-1}}$ and the same number $m$  of  ovals  from
$P^{\alpha_n}$. Moreover  each of these ovals contains exactly 2 points
of $P^{\alpha_{n-1}}\cap P^{\alpha_{n}}$. Thus $A=\cup_{w\in[W]} w(a)$
consists of  $m\cdot  2^{n-2}$  ovals  of  $\alpha_{n-1}$  and $m\cdot
2^{n-2}$ ovals of $\alpha_n$.

     Suppose $P\setminus A$ is connected.  Then we contract to  a
point each boundary component of $P\setminus A$. Thus  we obtain a compact
surface $P'$ of genus $g'=g(P')=g-(m+2)\cdot 2^{n-2}$. The forms
$\alpha_1,\dots,\alpha_n$ give involutions  $\alpha'_i:P'\rightarrow
P'$. Consider a complex structure on $P'$ such that
$\alpha'_1,\dots,\alpha'_n$ are real forms of $P'$. They generate a
commutative equipment such that
$$\sum^n_{i=1}\vert\alpha'_{i}\vert=\sum^n_{i=1}\vert\alpha_{i}\vert-m\cdot
2^{n-1}.$$ Therefore $$h(P,W)=\sum^n_{i=1}\vert\alpha_{i}\vert-2g=
\sum^n_{i=1}\vert\alpha'_{i}\vert+m\cdot 2^{n-1}-2g'-m\cdot 2^{n-1}-2^{n+1}
<h(n,g')\leqslant f(n).$$ If $P\setminus A$ is disconnected then it
forms 4 connected components. Let $P_1$ be one of them.
Contract to a point each boundary component of $P_1$ to obtain
a compact surface $P'$ with
$$g'=g(P')=\frac{1}{4}(g-(m-1)2^{n-2}- 4(2^{n-2}-1)).$$ The forms
$\alpha_1,\dots,\alpha_{n-2}$ give involutions
$\alpha'_i:P'\rightarrow P'$.  Consider a complex  structure  on  $P'$
such that $\alpha'_1,\dots,\alpha'_{n-2}$ are real forms of $P'$. They
generate a commutative equipment such that
$$\sum^{n-2}_{i=1}\vert\alpha'_{i}\vert=\frac{1}{4}(\sum^n_{i=1}
\vert\alpha_{i}\vert-m\cdot 2^{n-1}).$$
Thus $$h(P,W)=\sum^n_{i=1}\vert\alpha_{i}\vert-2g=
4\sum^{n-2}_{i=1}\vert
\alpha'_i\vert +m\cdot 2^{n-1}-8g'-(m-1)2^{n-1}-8(2^{n-2}-1)=$$
$$4(\sum^{n-2}_{i=1}\vert\alpha'_{i}\vert-2g')-3\cdot 2^{n-1}+8=
4f(n-2)-3\cdot 2^{n-1}+8.$$ If $n>4$ then
$$h(P,W)\leqslant 4(-(n-11)2^{n-5}-2)-3\cdot
2^{n-1}+8$$
$$=2^{n-3}(-n+11-12)-8+8=2^{n-3}(-n-1)<-2^{n-3}(n-9)-2=f(n).$$
If $n=4$  then $$h(P,W)\leqslant 4\cdot 2- 3\cdot 2^3+8=-8<f(4).$$ If
$n=3$ then
$$P^{\alpha_1}\cap (P^{\alpha_2}\cup   P^{\alpha_3})=\emptyset$$   and
lemma 3.3 follows from lemma 3.2. $\square$

\vskip 0.4cm
\noindent
{\bf Lemma 3.4.} {\sl Let $(P,W)$ be a commutative equipment such that
$[W]=(\alpha_1,\dots,\alpha_{n})$, where  $n>2$.  Then
$(n-4)2^{n-3}\leqslant g-1$. Moreover if every connected component of $P^W$
does not belong to $P^{\alpha_i}\cup P^{\alpha_j}$ for each $i,j$,  then
$$\sum_{i=1}^n\vert\alpha_i\vert\leqslant 2g(P)-(n-9)2^{n-3}-2.$$}

\noindent
{\sl Proof:} Let $C\subset P$ be a camera of $W$ and
$\widetilde a_1,\dots, \widetilde a_k$ be its boundary contours. The
contour $\widetilde a_i$ consists of segments
$$\ell_{i1},\dots,\ell_{im_i},$$
where $$\ell_{ij}\subset P^{\sigma(i,j)},\,\sigma(i,j)\in [W]$$
and $$\sigma(i,j)\ne \sigma(i,j+1),\,\sigma(i,1)\ne \sigma(i,m_i).$$
Let $t_i$ be the number of distinct  elements between
$$\sigma(i,1),\dots,\sigma(i,m_i).$$ Our conditions give $t_i\geqslant
3$. Moreover $$\sum^k_{i=1}t_i\geqslant n.$$ Put
$$\sigma(i,j+nm_i)=\sigma(i,j)$$ for $n\in Z$. Consider
$$L^1_i=\{\ell_{ij}\vert\sigma(i,j-1)=\sigma(i,j+1)\},$$
$$L^2_i=\{\ell_{ij}\vert\sigma(i,j-1)\ne\sigma(i,j+1)\}.$$
Let $s_i$ be the number of elements in $L^2_i$. It follows from
$t_i\geqslant 3$ that $s_i\geqslant t_i-1$.

Let $\widetilde P=P/W$ and $\varphi:P\rightarrow\widetilde  P$  be  the
natural projection.  Then  $\varphi^{-1}(\ell_{ij})$ consists of ovals
of $\sigma(i,j)$.  The  number  of  these  ovals   is   $2^{n-2}$   if
$\ell_{ij}\in L^1_i$ and $2^{n-3}$ if $\ell_{ij}\in L^2_i$. Thus
$$\sum^n_{i=1}\vert\alpha_i\vert=\sum^k_{i=1}(s_i\cdot 2^{n-3}+
(m_i-s_i)\cdot 2^{n-2})=\sum^k_{i=1}(m_i\cdot 2^{n-2}- s_i\cdot
2^{n-3})\leqslant$$
$$\sum^k_{i=1}(m_i\cdot 2^{n-2}-(t_i-1)2^{n-3})=(\sum^k_{i=1}m_i)\cdot
2^{n-2}+k\cdot 2^{n-3}-(\sum^k_{i=1} t_i)\cdot 2^{n-3}\leqslant$$
$$(\sum^k_{i=1}m_i)\cdot 2^{n-2}+(k-n)\cdot 2^{n-3}.$$

On the  other hand,  it follows from theorem 2.3 that
$$(P,W)=(\Lambda/\text{Ker}\ \psi,\ G/\text{Ker}\ \psi),$$ where
$(G,\varphi)$ is a $g-$planar realization the swelling Coxeter
system $(W,  [W],  T)$,  and $$T=\{\sigma(i,j), i=1,\dots,k,
j\in Z\}.$$ It follows from Riemann-Hurwitz's theorem ([9], 4.14.21)
that $$4g-4=2^n(4\widetilde   g-4+2k+\frac{1}{2}\sum^k_{i=1}   m_i),$$
where $g=g(P), \quad \widetilde g=g(\widetilde P)$. Thus
$$g-1\geqslant 2^{n-2}(-2+\frac{1}{2}n) =2^{n-3}(n-4)$$ and
$$\sum^n_{i=1}\vert\alpha_i\vert-2g\leqslant(\sum^k_{i=1}    m_i)\cdot
2^{n-2}+(k-n)\cdot 2^{n-3}   -2^{n-1}(4\widetilde   g-4+2k+\frac{1}{2}
\sum^k_{i=1} m_i)-2\leqslant$$
$$\leqslant  (k-n)\cdot  2^{n-3}+2^{n+1}-k\cdot
2^n-2\leqslant (1-n)\cdot 2^{n-3}+  2^{n+1}- 2^n-2=$$
$$=-2^{n-3}(n-1-16+8)-2=-(n-9)\cdot 2^{n-3}-2. \square$$

\vskip 0.4cm
\noindent
{\bf Lemma 3.5.} {\sl Let $(P,W)$ be a real orientable equipment and
$\tau_1,\dots,\tau_n\in [W]$ (where $n>2$) be  non-conjugate. Then
there exists a commutative equipment $W'\subset W$ such that
$$[W']=(\alpha_1,\dots,\alpha_n),\quad \alpha_i\ne\alpha_j$$ and
$$\sum^n_{i=1}\vert\tau_i\vert\leqslant \sum^n_{i=1}\vert
\alpha_i\vert.$$}

\noindent
{\sl Proof:} Let $\widetilde W$ be the real orientable equipment
generated by $\tau_1,\dots,\tau_n$ and $C$ its camera.  It follows from
lemma 1.1 that there exist  $\beta_1,\dots,\beta_n\in  S_C$  such  that
$\beta_i=w_i\tau_i w^{-1}_i$, $w_i\in \widetilde W$. Let $\widetilde W'$
be the  real  orientable  equipment,  generated by  $\widetilde
S'=(\beta_1,\dots,\beta_n)$. Put  $\alpha_j=\beta_j$  if  $\beta_j$
belongs to the center of $\widetilde W'$.

Let us now assume that $\beta_j$ does not belong to the center of
$\widetilde W'$.  It  follows  from  theorem 1.1 that $(\widetilde W',
\widetilde S')$ is a Coxeter  system.  Moreover  $\beta_1,\dots,
\beta_n$ are non-conjugated in $\widetilde W'$.
We observe, using the
classification of Coxeter systems [2, VI, $\S$4], that
there exists only one $\beta_i\in \widetilde S'$
such that $\beta_i\beta_j\ne\beta_j\beta_i$. For similar reasons
$\beta_i\beta_k=\beta_k\beta_i$ if $k\ne j$.
The order $2m$ of $\beta_i\beta_j$ is even because $\beta_i$ and $\beta_j$ are
non-conjugated in $\widetilde W'$.
Put $\gamma=(\beta_i\beta_j)^m$.  Then $\gamma$
belongs to the center of $\widetilde W'$. For $\vert\beta_i\vert\geqslant\vert
\beta_j\vert$ we put  $\alpha_i=\beta_i$,  $\alpha_j=\gamma\beta_i$.  For
$\vert\beta_j\vert\geqslant\vert\beta_i\vert$ we put  $\alpha_j=\beta_j$,
$\alpha_i=\gamma\beta_j$. Then  $\alpha_1,\dots,\alpha_n$  generate  a
commutative equipment  $W'$ and $$\sum^n_{i=1}\vert\alpha_i\vert
\geqslant \sum^n_{i=1}\vert\tau_i\vert. \square$$

\vskip 0.4cm
\noindent
{\bf Theorem 3.1.} {\sl Let $(P,W)$ be a real orientable equipment, $g(P)=g$ and
$\tau_1,\dots,\tau_n\in [W]$  (where $n>2$) be non-conjugated in
$W$. Then  $(n-4)2^{n-3}\leqslant g-1$ and
$$\sum^n_{i=1}\vert\tau_i\vert\leqslant 2g-(n-9)2^{n-3} -2.$$}

\noindent
{\sl {Proof:}} Due to lemma 3.5 it suffices to
prove theorem 3.1 for commutative equipment.  For this case it follows
from lemma 3.4 that $(n-4)2^{n-3}\leqslant g-1$.
We use an induction on $g=g(P)$ for to prove
$$\sum^n_{i=1}\vert\tau_i\vert\leqslant 2g-(n-9)2^{n-3} -2.$$
If $g(P)=0$ then $n=3$ and $\sum^n_{i=1}\vert\tau_i\vert=3< 6-2$.
Let us assume that the statement
is proved for the cases $g(P)< g$. If it exists a connected
component $a$ of $P^W$ such that $a\subset P^{\tau_i}\cup  P^{\tau_j}$,
then the statement of theorem 3.1 follows from lemma 3.3. Otherwise it
follows from lemma 3.4. $\square$

\vskip 0.4cm
\noindent
{\bf Theorem 3.2.} {\sl For any $n>2$ and $m\geqslant 0$, 
where $n+2m>4$, there exists a
commutative equipment  $(P,W)$ such that the forms
$[W]=(\tau_1,\dots,\tau_n)$ are non-conjugated with
respect to  holomorphic automorphisms  of $P$,
$$g(P)=2^{n-3}(n+2m-4)+1$$ and
$$\sum^n_{i=1}\vert\tau_i\vert=2g(P)-(n-9)2^{n-3}-2.$$}

\noindent
{\sl Proof:} Consider a rectangular $(n+2m)$-gon
with incommensurable sides  $\ell_1,\dots,\ell_{n+2m}$
on Lobachevskij plane $\Lambda$. Let
$G\subset \overline{\text{Aut}}
 (\Lambda)$ be  the group generated by
reflections $\sigma_i$ in the sides $\ell_i$ of the polygon. Let
$\widetilde W\cong(Z_2)^n$ be the group generated by $n$ involutions $s_1,\dots,s_n$.
Consider the epimorphism $\psi:G\rightarrow \widetilde W$ such that
$\psi(\sigma_i)=s_i$ for  $i=1,\dots,n-1$,  $\psi(\sigma_{n+2j})=s_2$,
$\psi(\sigma_{n+2j+1})=s_n$ for $j=0,...,m-1$.  Then $(G,\psi)$  is  a
plan realization. It follows from lemma 2.1 that $(P,W)=(\Lambda/\text{Ker}
\ \psi, \quad G/ \text{Ker}\ \psi)$ is a commutative equipment and
$[W]=(\tau,\dots,\tau_n)$, where
$$\tau_i=\sigma_i / \text{Ker}\ \psi\quad\text{for}\ i=1,\dots,n-1, \quad
\tau_2=\sigma_{n+2i} / \text{Ker}\ \psi \quad \text{for}\ i=0,\dots, m,$$
$$\tau_n=\sigma_{n+2j+1}/\text{Ker}\ \psi,
\quad  \text{for}\  j=0,\dots,m-1.$$   Riemann-Hurwitz's
formula [9, 4.14.11] gives $$g(P)=2^{n-3}(n+2m-4)+1.$$
Let  $\psi:P\rightarrow  P/W$  be the natural projection.
Then $\psi^{-1}(\ell_i)$ forms $2^{n-3}$ ovals of $\tau\in [W]$ for
$i=2,\dots,n$ and it forms $2^{n-2}$ ovals of $\tau\in [W]$ for $i=1$ and $i>n$.
Thus $$\sum^n_{i=1}\vert\tau_i\vert=(n-1)2^{n-3}+(2m+1)2^{n-2}$$ and
$$2g(P)-\sum^n_{i=1}\vert\tau_i\vert=2^{n-3}(2n+4m-8)+2-(n-1)2^{n-3}-(4m+
2)2^{n-3}=$$ $$=(n-9)2^{n-3}+2.$$ It follows from incommensurability  of
$\ell_i$ that $\tau_i$ are  non-conjugated with respect to
holomorphic automorphisms (i.e.  isometries with respect to
Lobachevskij metric) of $P$. $\square$

\vskip 0.4cm
\noindent
{\bf Corollary 3.1.} {\sl Let $X_1,\dots,X_n \ (n>3)$ be orientable
non-singular, non-reducible real algebraic  curves
of genus $g>1$ such that for any $i\ne j$ $X_i$ is
non-isomorphic to $X_j$ over $\Bbb R$ but isomorphic over
$\Bbb C$. Then
$$\sum^n_{i=1}\vert X_i\vert\leqslant 2g-(n-9)2^{n-3}-2 \quad\text{and}\quad
(n-4)2^{n-3}\leqslant g-1$$
and this estimate is attained  for each $n$ for infinite
number of $g$.}

\noindent
{\sl Proof:} By definition there exists a Riemann surface $P$ and
biholomorphic maps   $\psi:P_i\rightarrow P$, such that
$X_i=(P_i,\alpha_i)$. Then $\tau_i=\psi_i\alpha_i\psi^{-1}_i$ generate
an orientable equipment of $P$.  Thus corollary 3.1 follows from
theorems 3.1 and 3.2. $\square$

\centerline {\bf References}

1. Alling N.L., Greenleaf N.,  Foundations of the theory of Klein
surfaces, {\sl {Lecture Notes in Math.}}, {\bf 219} (1971).

2. Bourbaki N., {\sl {Groupes et Alg\'ebres de Lie, Chapitres 4,5 et 6}},
Hermann, Paris (1968).

3. Bujalance E., Gromadzki G., Singerman D., On the number of real curves
associated to a complex algebraic curve, {\sl {Proc.  amer.  math.  soc.}}
{\bf 120} (2) (1994), 507-513.

4. Natanzon S.M., Finite groups of homeomorphisms of surfaces ans real
forms of complex algebraic curves, {\sl {Trans. Moscow Math. Soc.}}
{\bf 51} (1989), 1-51.

5. Davis M., Groupes generated by reflections and aspherical
manifolds not covered by Euclidean space,
{\sl {Ann. Math.}} {\bf 117} (1983), 293-324.

6. Natanzon S.M., Real borderings of complex algebraic curves and
Coxeter group, {\sl {Russian Mat. Surveys}} {\bf 51} (6) (1996),
1216-1217.

7. Natanzon S.M.,  A Harnack-type theorem for a family of
complex-isomorphic real algeb\-raic curves,
{\sl {Russian Math. Surveys}}, {\bf 52} (6) (1997).

8. Singerman D., Mirrors on Riemann surfaces. {\sl {Contemporary
Mathematics}} {\bf 184} (1995), 411-417.

9. Zieschang H., Vogt E., Coldewey H.-D., 1980,  Surfaces and planar
discontinuous groups, {\sl {Lecture Notes in Math.}} {\bf 835} (1980).

\vskip 0.4cm

e-mail: natanzon\@mccmc.ru

\end